\documentclass{amsart}
\usepackage{amscd,amssymb}
\usepackage[all]{xy}

\newtheorem{theorem}{Theorem}[section]
\newtheorem{prop}[theorem]{Proposition}
\newtheorem{sps}[theorem]{Supposition}
\newtheorem{corol}[theorem]{Corollary}
\newtheorem{conj}[theorem]{Conjecture}
\theoremstyle{definition}
\newtheorem{rmk}[theorem]{Remark}
\newtheorem{defin}[theorem]{Definition}
\theoremstyle{remark}
\newtheorem*{erem}{Remark}
\numberwithin{equation}{section}

\def\iff{if and only if }

\def\bul{_\bullet}	\def\vR{\mathbf r_\bullet}
\def\mdl#1{\mathrm{mod}\mbox{-}#1}
\def\Mdl#1{\mathrm{Mod}\mbox{-}#1}
\def\Mk{\Bbbk}	\def\8{\infty}
\def\={\setminus}	\def\*{\otimes}
\def\bop{\bigoplus}  \def\+{\oplus}
\def\larr{\longrightarrow}  \def\ol{\overline}
\def\lst#1#2{ #1_1 , #1_2 , \dots , #1_{#2} }
\def\row#1#2{( #1_1 , #1_2 , \dots , #1_{#2}) }
\def\xx{\times}

\def\hom{\mathop{\mathcal H\!\mathit{om}}\nolimits}
	\def\CS{coherent sheaf}
\def\pam{\mathop\mathrm{par}\nolimits}

\def\set#1{\left\{\,#1\,\right\}}
\def\setsuch#1#2{\left\{\,#1\,|\,#2\,\right\}}
\def\dlim{\varinjlim}	 \def\ilim{\varprojlim}
\def\im{\mathop{\mathrm{Im}}\nolimits}
\def\Ker{\mathop{\mathrm{Ker}}\nolimits}
\def\Hom{\mathop{\mathrm{Hom}}\nolimits}
\def\rad{\mathop{\mathrm{rad}}\nolimits}
\def\Aut{\mathop{\mathrm{Aut}}\nolimits}
\def\rk{\mathop{\mathrm{rk}}\nolimits}

\def\kA{\mathcal A} \def\kB{\mathcal B}
\def\kC{\mathcal C}  \def\kH{\mathcal H}
\def\kK{\mathcal K}  \def\kP{\mathcal P}
\def\kV{\mathcal V}	  \def\kR{\mathcal R}
\def\kG{\mathcal G}	 \def\kI{\mathcal I}	  
\def\kO{\mathcal O}	\def\kL{\mathcal L}
\def\kJ{\mathcal J}	\def\kE{\mathcal E}	
\def\bA{\mathbf A}	\def\bB{\mathbf B}
\def\bR{\mathbf R}	\def\bI{\mathbf I}
\def\bJ{\mathbf J}	\def\bK{\mathbf K}
\def\fB{\mathbf b}	\def\fA{\mathbf a}
\def\fR{\mathbf r}	\def\dM{\mathfrak M}
\def\mD{\mathbb D} \def\mP{\mathbb P}
\def\om{\omega}	\def\be{\beta}
\def\al{\alpha}	\def\la{\lambda}
\def\ga{\gamma}

\begin{document}
\title{Semi-continuity for derived categories}
\author{Yuriy Drozd}
\address{Kyiv Taras Shevchenko University, Department of Mechanics and Mathematics,
 01033 Kyiv, Ukraine}
\email{yuriy@drozd.org}
\thanks{ This work was inspired by my visit to the University of California at Santa Barbara and
 accomplished during my stay at the University of Kaiserslautern. I thank the host Universities and the
 supporting foundations, CRDF (Award UM2-2094) and DFG (Schwerpunkt ``Globale Methoden
 in der komplexen Geometrie'') for this opportunity. I am also grateful to 
 Birge Huisgen-Zimmermann and Igor Burban for useful discussions.}
\begin{abstract}
 We prove that the number of parameters defining a complex of projective modules over an
 algebra is upper semi-continuous in families of algebras. Supposing that every algebra is either derived
 tame or derived wild, we get that a degeneration of a derived wild algebra is also derived wild. The proof
 follows the pattern of the paper \cite{dg} and rests upon universal families with projective bases. We also
 explain why the so-called counter-example of Br\"ustle \cite{br} is in fact not a counter-example.
\end{abstract}
\maketitle


 In the representation theory of finite dimensional algebras it is usual to distinguish three types of algebras:
 representation finite, tame and wild (cf. \cite{gr,d2}). A useful tool in establishing representation type is
 provided by deformation theory. It is based on upper semi-continuity of \emph{parameter number}
 $\pam(n,\bA(x))$ defining an $n$-dimensional module over an algebra
 $\bA(x)$, when $\bA(x)$ is an algebraic family of algebras (see \cite{gd,dg}). During last years analogous
 investigation of derived categories has been started, especially derived tame and wild algebras have been considered.
 In this article we define parameter numbers of complexes and prove that they are also upper semi-continuous
 in families of algebras.  We follow the technique elaborated in \cite{dg} and going back to a paper of
 H.\,Kn\"orrer \cite{kn}. Especially, our proof depends on the construction of 
 (almost) universal families of complexes with projective bases. As a corollary, we prove that if an algebra,
 which is not derived tame, degenerates to another algebra, the latter is also not derived tame
 (note that most people working on the subject believe that `not tame' means `wild' in this situation too).
 Recently Th.\,Br\"ustle \cite{br} has announced a counter-example to the last assertion. We explain why it is actually
 \emph{not} a counter-example.

\section{Categories $\kK^n(\bA)$}
 \label{s1}

 Let $\bA$ be a ring. We denote by $\Mdl\bA$ ($\mdl\bA$) the category of right $\bA$-modules
 (respectively, of finitely generated $\bA$-modules).
 We define the category $\kK^n(\bA)$ as follows.
\begin{enumerate}
\item 
Its objects are finite complexes of  projective $\bA$-modules 
 \begin{equation}\label{e11}
 \begin{CD}
 P\bul:\quad P_n@>{d_n}>>P_{n-1}@>{d_{n-1}}>>\dots\larr P_m   
 \end{CD} 
 \end{equation}
 ($m\le n$). We set $P_k=0$ for $k<m$.
\item 
Morphisms in $\kK^n(\bA)$ are homomorphisms of complexes modulo
\emph{quasi-homo\-topy}. Namely, two homomorphisms $\phi=(\phi_k)$ and
$\psi=(\psi_k)$ from $P\bul$ to $P'\bul$ are called
\emph{quasi-homotopic} if there are homomorphisms of modules
$s_k:P_k\to P'_{k+1}$ such that
$\psi_k=\phi_k+s_{k-1}d_k+d'_{k+1}s_k$ for all $k< n$. We write $\phi\stackrel n\sim\psi$
 in this case and call $s=(s_k)$ a \emph{quasi-homotopy from $\phi$ to $\psi$}.
\end{enumerate}
 Note that the last homomorphisms $\phi_n,\psi_n$ do not influence quasi-homotopy at all.

There is a natural functor $\kI_n:\kK^n(\bA)\to \kK^{n+1}(\bA)$. Namely, if $P\bul$ is a complex
 from $\kK^n$, choose an epimorphism $d_{n+1}:P_{n+1}\to P_n$ with $\im d_{n+1}=\Ker d_n$
 and define $\kI_n P\bul$ as the complex 
 $$
 \begin{CD}
     P_{n+1}@>{d_{n+1}}>> P_n @>{d_n}>> \dots\larr  P_m.
 \end{CD} 
 $$ 
 If $\phi$ is a homomorphism $P\bul\to P'\bul$, we can lift it to a homomorphism
 $\kI_n\phi:\kI_n P\bul\to \kI_n P'\bul$ in a usual way. Moreover, if $\phi\stackrel n\sim\psi$ and $s=(s_n)$
 is a corresponding quasi-homotopy, one easily sees that $d'_n(\phi_n-\psi_n-s_{n-1}d_n)=0$, hence there is
 a mapping $s_n:P_n\to P'_{n+1}$ such that $\phi_n-\psi_n-s_{n-1}d_n=d'_{n+1}s_n$, so
 we get a quasi-homotopy $\kI_n\phi\stackrel{n+1}\sim\kI_n\psi$. Therefore the functor $\kI_n$ is well
 defined. It gives rise to the direct limit $\kK^\om(\bA)=\dlim_n\kK^n(\bA)$. 

 \begin{prop}\label{11}
  The category $\kK^\om(\bA)$ is equivalent to the bounded derived category $D^b(\Mdl\bA)$.
  If $\bA$ is noetherian, the bounded derived category $D^b(\mdl\bA)$ is equivalent to the full
 subcategory $\kK^\om_f(\bA)=\dlim_n\kK^n_f(\bA)$ of $\kK^\om(\bA)$, where $\kK^n_f(\bA)$ is
 the full subcategory of $\kK^n(\bA)$ consisting of complexes of finitely generated modules. 
 \end{prop} 
 (This result traces back to the paper \cite{bs}.) 
\begin{proof}
   Consider the functor $\kJ_n:\kK^n(\bA)\to D^b(\Mdl\bA)$, which maps a complex $P\bul$ to the
 complex
 $$
  \begin{CD}
 \Ker d_n\larr P_n@>{d_n}>>P_{n-1}@>{d_{n-1}}>> \dots\larr P_m.  
 \end{CD} 
  $$ 
 One can check that $\kJ_{n+1}\kI_n\simeq\kJ_n$, so we get the limit functor $\kJ=\dlim_n\kJ_n:
 \kK^\om(\bA)\to D^b(\Mdl\bA)$. On the other hand, let $C\bul$ be a complex from $D^b(\Mdl\bA)$
 such that $C_k=0$ for $k\ge n$. Consider its projective resolution $P\bul$; it is exact at all $P_k$
 with $k\ge n$. Let $P^{(n)}\bul\in \kK^n(\bA)$ be the complex that coincide with $P\bul$ for $k\le n$.
 Then $C\bul\simeq\kJ_nP^{(n)}\bul $ in $D^b(\Mdl\bA)$. Moreover, $\kI_n P^{(n)}\bul\simeq
 P\bul^{(n+1)}$, so $C\bul\simeq\kJ P^{(n)}$ and we get a functor $D^b(\Mdl\bA)\to \kK^\om(\bA)$
 inverse to $\kJ$.

 The assertion about noetherian case is obvious.
 \end{proof} 

 Note that there are also natural functors $\kE_n:\kK^{n+1}(\bA)\to\kK^n(\bA)$: we just omit the 
 term $P_{n+1}$. Hence the inverse limit $\kK^\8(\bA)=\ilim_n\kK^n(\bA)$ is defined and 
 the following result hold.

 \begin{prop}\label{12}
  The category $\kK^\8(\bA)$ is equivalent to the right bounded derived category $D^-(\Mdl\bA)$.
   If $\bA$ is noetherian, the right bounded derived category $D^-(\mdl\bA)$ is equivalent to the full
 subcategory $\kK^\8_f(\bA)=\ilim_n\kK^n_f(\bA)$ of $\kK^\8(\bA)$.
 \end{prop} 
 The proof is quite analogous to that of Proposition \ref{11} and we omit it.

\medskip
 Suppose now that the ring $\bA$ is semiperfect \cite{bass}; set $\bR=\rad\bA$.
 Then any right bounded complex of finitely generated projective $\bA$-modules is
 homotopic to a \emph{minimal complex} $P\bul$, i.e. such that
 $\im d_n\subseteq P_{n-1}\bR$ for all $n$. Thus, considering $\kK^n_f(\bA)$, we may confine
 ourselves to minimal complexes, and we shall always do so. Note that two minimal complexes
 are homotopic \iff they are isomorphic. We call a minimal complex  \eqref{e11} \emph{reduced}
 if $\Ker d_n\subseteq P_n$. Any complex from $\kK^n_f(\bA)$  is isomorphic (as complex) to
 a direct sum of a reduced one and a complex having all zero components except maybe the $n$th.
 The latter is a zero object of $K^n(\bA)$, so any complex is isomorphic in $\kK^n_f(\bA)$ to a reduced
 one. Moreover, if a homomorphism $\phi:P\bul\to P'\bul$ of reduced complexes from $\kK^n_f(\bA)$
 is quasi-homotopic to 0, $\im\phi_k\subseteq P'_k\bR$ for all $k$ including $k=n$. Therefore if two reduced
 complexes are isomorphic in $\kK^n_f(\bA)$, they are isomorphic as complexes. 

 We denote by $\kC^n(\bA)$ the category of minimal complexes $P\bul$ with $P_k=0$ for $k>n$ and
 by $\kC^n_0(\bA)$ its full subcategory of reduced complexes. Then the natural functor $\kC^n(\bA)\to
 \kK^n_f(\bA)$ induces a \emph{representation equivalence} of $\kC^n_0(\bA)$ onto $\kK^n_f(\bA)$.
 Here we call a functor $F:\kA\to\kB$ a \emph{representation equivalence}, if
 \begin{itemize}
\item   $F$ is \emph{dense}, i.e. every object from $\kB$ is isomorphic to $FA$ for some object $A\in\kA$;
 \item
  $F$ is \emph{full}, i.e. all induced mappings $\kA(A,A')\to\kB(FA,FA')$ are surjective;
 \item
  $F$ is \emph{conservative}, i.e. $F\phi$ is an isomorphism \iff $\phi$ is an isomorphism.
 \end{itemize}
 An evident consequence of these conditions is
 \begin{itemize}
 \item
  $FA$ is indecomposable \iff $A$ is indecomposable.
 \end{itemize}
 
 Let $\lst At$ be all pairwise non-isomorphic indecomposable projective $\bA$-modules (all of them are
 direct summands of $\bA$). If $P$ is a finitely generated projective $\bA$-module, it uniquely decomposes as
 $\bop_{i=1}^tp_iA_i$. We denote by $\fR(P)$ the vector $\row pt$ and call it the \emph{rank} of $P$.
 For any complex \eqref{e11} from $\kC^n(\bA)$ we denote by $\vR(P\bul)$
 and call the \emph{vector rank} of $P\bul$ the sequence $(\fR(P_n),\fR(P_{n-1}), \dots,\fR(P_m))$.
 As we have already remarked, every complex from $\kC^n(\bA)$ decomposes as
 $P\bul\+(\bop_{i=1}^ta_iA_i[n])$, where $P\bul$ is a reduced complex and $A[n]$ denotes, as
 usually, the complex with a unique non-zero component, namely the $n$th one, equal $A$. Thus, from the
 viewpoint of classification problem, there is no essential difference between $\kC^n(\bA)$ and $\kK^n_f(\bA)$.
 
 Given a vector $\fR=\row pt$, we denote $\fR A=\bop_{i=1}^tp_iA_i$ and set $\bA(\fR,\fR')=
 \Hom_\bA(\fR A,\fR'A),\ \bR(\fR,\fR')=\Hom_\bA(\fR A,\fR'A\bR)$. For any sequence
 $(\fR_n,\fR_{n-1},\dots,\fR_m)$ we consider the set $\kC^n(\fR)$ of minimal complexes
 \begin{equation}\label{e12}
  \begin{CD}
  \fR_nA@>{d_n}>>\fR_{n-1}A@>{d_{n-1}}>>\dots\larr \fR_mA,  
 \end{CD} 
 \end{equation} 
 or, the same, of sequences $(d_n,d_{n-1},\dots,d_{m+1})$ with $d_k\in\bR(\fR_k,\fR_{k-1})$
 such that $d_kd_{k+1}=0$ for all $m<k<n$. Two sequences $(d_k)$ and $(d'_k)$ define
 isomorphic complexes \iff there are invertible mappings $\al_k\in\bA(\fR_k,\fR_k)\ (m\le k\le n)$
 such that $\al_{k-1}d_k=d'_k\al_k$ for all $m<k\le n$. Especially two sequences $(d_n)$ and
 $(\la_n d_n)$, where $\la_n$ are invertible elements from the centre of $\bA$, always define isomorphic
 complexes.

 If $\bA$ is a finite dimensional algebra over a field $\Mk$, it allows us to consider complexes from
 $\kC^n(\bA)$ of a fixed vector rank $\vR$ as ($\Mk$-valued) points of an algebraic variety $\kC(\vR)$,
 which is a subvariety of $\kH(\fR\bul)=\prod_{k=m+1}^n\bR(\fR_k,\fR_{k+1})$. Moreover, 
 homothetic sequences $(d_n)$ and $(\la d_n)$ define isomorphic complexes, so considering the classification
 problem we may replace the vector space  $\kH(\fR\bul)$ by the projective space
 $\mP(\fR\bul)=\mP(\kH(\fR\bul))$ and $\kC(\fR\bul)$ by its
 image $\mD(\fR\bul)$ in $\mP(\fR\bul)$, which is a projective variety. 
 
 \section{Families of complexes}
 \label{s2}

 From now on let $\bA$ be a finite dimensional algebra over an algebraically closed field $\Mk$
 and $\bI\subseteq\bR$ be an ideal. 
 We define an \emph{$\bI$-family of $\bA$-complexes} over an algebraic variety $X$
 as a complex of flat \CS s of $\bA_X$-modules
 \begin{equation}\label{e21}
 \begin{CD}
   \kP\bul:\quad \kP_n@>{d_n}>>\kP_{n-1}@>{d_{n-1}}>>\dots\larr \kP_m ,  
 \end{CD} 
 \end{equation} 
 where $\bA_X=\kO_X\*\bA$, such that $\im d_k\subseteq\kP_{k-1}\bI$ for all $m< k\le n$.
 If $\bI=\bR$, we also call  such a complex a \emph{family of minimal $\bA$-complexes}.
  Given a family $\kP\bul$ and a point $x\in X$, we get a complex
 $\kP\bul(x)$ from $\kC^n(\bA)$. Since $\kP_n$ are locally free over $\kO_X$, the ranks
 $\fR(\kP_k(x))$ are locally constant. We usually suppose $X$ connected; then these ranks
 are constant, so we can define $\fR(\kP_k)$ and $\fR\bul(\kP\bul)$. 
 Consider the set
 $$
  \kI=\setsuch{(x,y)\in X\xx X}{\kP\bul(x)\simeq\kP\bul(y)}.
 $$
 It is a constructible subset of $X\xx X$ and for each $x\in X$ the set
 $$
  \kI(x)=\pi_1^{-1}(x)\cap\kI\simeq\setsuch{y\in X}{\kP\bul(y)\simeq\kP\bul(x)}
 $$
 is closed in $\kI$, hence also constructible. One can easily derive
 from the standard results on dimensions of fibres \cite[Exercise II.3.22]{ha} that the set
 $X_i(\kP\bul)=\setsuch{x\in X}{\dim\kI(x)\le i}$ is constructible too. (It is also a consequence
 of Propositions \ref{21} and \ref{23} below.) We define the
 \emph{number of parameters} in the family $\kP\bul$ as 
 $$ 
   \pam(\kP\bul)=\max_i\set{\dim X_i(\kP\bul)-i},
 $$ 
 the \emph{number of parameters in $\bI$-families of vector rank} $\fR\bul$ as 
 $$ 
   \pam(\fR\bul,\bA,\bI)=\max\setsuch{\pam(\kP\bul)}{\fR\bul(\kP\bul)=\fR\bul},
 $$ 
 and the \emph{number of parameters in families of minimal complexes} of vector rank $\fR$ as 
 $$ 
   \pam(\fR\bul,\bA)=\pam(\fR\bul,\bA,\bR).
 $$ 
 It is a formal version of the intuitive impression about the number of parameters necessary to define an
 individual complex inside the family. Of course, we are mainly interested in the ``absolute'' value
 $\pam(\fR,\bA)$, but further on we need also its ``relative'' version. Obviously 
 $\pam(\fR\bul,\bA,\bJ)\le\pam(\fR\bul,\bA,\bI)$ if $\bI\supseteq\bJ$, especially always
 $\pam(\fR\bul,\bA,\bI)\le\pam(\fR\bul,\bA)$.
 
 A family of complexes \eqref{e21} is called \emph{non-degenerate} if for every point $x\in X$ at least one
 of the homomorphisms $d_k(x)$ is non-zero. Obviously, there is an open set $U\subseteq X$ such that the
 restriction of $\kP\bul$ onto $U$ is non-degenerate, and considering classification problems, as well as
 calculating parameter numbers, it is enough to deal with non-degenerate families. 

 We are able, just as in \cite{dg}, to construct some ``almost universal'' non-degenerate families.
 It is important that their bases are \emph{projective varieties}.  Namely, fix a vector rank $\fR\bul$ and
 set $\kH=\kH(\fR,\bI)=\bop_{k=m+1}^n\bI(\fR_k,\fR_{k-1})$, where
 $\bI(\fR,\fR')=\Hom_\bA(\fR A,\fR'A\bI)$. Consider the projective space $\mP=\mP(\fR\bul,\bI)=
 \mP(\kH)$ and its closed subset $\mD=\mD(\fR\bul,\bI)\subseteq\mP$ consisting of  sequences $(h_k)$
 such that $h_{k+1}h_k=0$ for all $k$. Because of the universal property of projective spaces
 \cite[Theorem II.7.1]{ha}, the embedding $\mD(\fR\bul) \to\mP(\fR\bul)$ gives rise to a non-degenerate
 $\bI$-family $\kV\bul=\kV\bul(\fR\bul,\bI)$ 
 \begin{equation*}\label{e22} 
 \begin{CD} 
    \kV\bul(\fR\bul):\quad \kV_n@>{d_n}>>\kV_{n-1}@>{d_{n-1}}>>\dots\larr \kV_m ,
 \end{CD} 
 \end{equation*} 
 where $\kV_k=\kO_\mD(n-k)\*\fR_kA$ for all $m\le k\le n$.
 We call $\kV\bul(\fR\bul,\bI)$ the \emph{canonical $\bI$-family of $\bA$-complexes} over $\mD(\fR\bul,\bI)$.
 Moreover, morphisms $\phi:X\to\mD(\fR\bul,\bI)$ correspond to non-degenerate families of shape
 \eqref{e21} with $\kP_k=\kL^{\*(n-k)}\*\fR_kA$ for some invertible sheaf $\kL$
 over $X$. Namely, such a family can be obtained as $\phi^*\kV\bul$ for a uniquely defined morphism
 $\phi$. 

  \begin{prop}\label{21}
  For every non-degenerate family of $\bI$-complexes $\kP\bul$ of vector rank $\fR\bul$ over an algebraic
 variety $X$ there is a finite covering open $X=\bigcup_jU_j$ such that the restriction of $\kP\bul$ onto each $U_j$
 is isomorphic to $\phi_j^*\kV(\fR\bul,\bI)$ for some morphism $\phi_j:U_j\to \mD(\fR\bul,\bI)$.
 \end{prop}
  \begin{proof}
 For each $x\in X$   there is an open neighbourhood $U\ni x$ such that all restrictions $\kP_k|U$ are isomorphic to
 $\kO_U\*\fR_kA$, so the restriction $\kP\bul|U$ is obtained from a morphism $U\to\mD(\fR\bul,\bI)$. 
 Evidently it implies the assertion.
 \end{proof} 
 Note that morphisms $\phi_j$ are not canonical, so we cannot glue them into a global morphism $X\to\mD(\fR\bul,\bI)$.

   \begin{corol}\label{22}
  $\pam(\fR\bul,\bA,\bI)=\pam(\kV\bul(\fR\bul,\bI))$.
 \end{corol} 

 The main advantage of the families $\kV\bul(\fR\bul,\bI)$ is the following.

 \begin{prop}\label{23}
  All sets
 $$ 
   \mD_i(\fR\bul,\bI)=\setsuch{x\in\mD(\fR\bul,\bI)}{\dim\setsuch{y\in\mD(\fR\bul,\bI)}
  {\kV\bul(\fR\bul,\bI)(y)\simeq\kV\bul(\fR\bul,\bI)(x)}\le i}
 $$ 
 are closed in $\mD(\fR\bul,\bI)$.
 \end{prop} 
 \begin{proof}
   Consider the group $G=G(\fR\bul)=\prod_{k=m}^n \Aut_\bA(\fR_kA)$. It acts on $\kH(\fR\bul,\bI)$:
 $(g_k)(h_k)=(g_{k-1}h_kg_k^{-1})$, hence also on $\mP(\fR\bul,\bI)$ and on $\mD(\fR\bul,\bI)$. Moreover,
 complexes $\kV(\fR\bul,\bI)(x)$ and $\kV(\fR\bul,\bI)(y)$ are isomorphic \iff the points $x$ and $y$ are in the same
 orbit of the group $G$. Hence $\mD_i(\fR\bul,\bI)=\setsuch{x\in\mD(\fR\bul,\bI)}{\dim Gx\le i}$, and it is well
 known that this set is closed.
 \end{proof} 

 In the next section we shall mainly consider complexes and families of  free modules, so we introduce corresponding
 notations. Let $\fA=\row at=\fR(\bA)$. For any sequence of integers $\fB=(b_n,b_{n-1},\dots,b_m)$ we set
 $\fB\fA=(b_n\fA,\dots,b_m\fA)$, $\mD(\fB,\bA,\bI)=\mD(\fB\fA,\bI)$ and $\pam(\fB,\bA,\bI)=
 \pam(\fB\fA,\bA,\bI)$, in particular $\pam(\fB,\bA)=\pam(\fB\fA,\bA,\bR)$.
 For any  $\fR=\row rt$ we denote by $]\fR/\fA[$ the smallest integer $b$ such that $ba_i\ge r_i$ for all $i$. If 
 $\fR\bul=(\fR_n,\dots,\fR_m)$, we set $]\fR\bul/\fA[=(]\fR_n/\fA[,\dots,]\fR_m/\fA[)$. Just analogous
 the values $[\fR/\fA]$ and $[\fR\bul/\fA]$ are defined. If $\fB=]\fR\bul/\fA[$ and $\fB'=[\fR\bul/\fA]$, then
 evidently 
 $$
  \pam(\fB'\fA,\bA,\bI)\le \pam(\fR\bul,\bA,\bI)\le \pam(\fB\fA,\bA,\bI).
 $$
 Therefore, when we are interested in the \emph{asymptotic} of the function $\pam(\fR\bul,\fA)$ for big ranks,
 we may only consider complexes of free $\bA$-modules.

 \section{Families of algebras}

  A \emph{family of algebras} over an algebraic variety $X$ is a sheaf $\kA$ of $\kO_X$-algebras, which is
 coherent and flat (thus locally free) as a sheaf of $\kO_X$-modules. For such a family and every sequence
 $\fB=(b_n,b_{n-1},\dots,b_m)$ one can define the function $\pam(\fB,\kA,x)=\pam(\fB,\kA(x))$.
 Our main result is the upper semi-continuity of these functions.
 
  \begin{theorem}\label{31}
  Let $\kA$ be a family of algebras based on a variety $X$. The set $X_j=\setsuch{x\in X}{\pam(\fB,\kA,x)\ge j}$
 is closed for every $\fB$ and every integer $j$.
 \end{theorem} 
 \begin{proof}
  We may assume that $X$ is irreducible. Let $\bK$ be the field of rational functions on $X$. We consider it as a
 constant sheaf on $X$. Set $\bR=\rad(\kA\*_{\kO_X}\bK)$ and $\kR=\bR\cap\kA$. It is a sheaf of nilpotent
 ideals. Moreover, if $\xi$ is the generic point of $X$, the factor algebra $\kA(\xi)/\kR(\xi)$ is semisimple. Hence
 there is an open set $U\subseteq X$ such that $\kA(x)/\kR(x)$ is semisimple, thus $\kR(x)=\rad\kA(x)$ for every $x\in U$.
 Therefore $\pam(\fB,\kA,x)=\pam(\fB,\kA(x),\kR(x))$ for $x\in U$, so $X_j=X_j(\kR)\cup X'_j$, where
 $$ 
   X_j(\kR)=\setsuch{x\in X}{\pam(\fB,\kA(x),\kR(x))\ge j}
 $$ 
 and $X'=X\=U$ is a closed subset in $X$. Using noetherian induction, we may suppose that $X_j'$ is closed,
 so we only have to prove that $X_j(\kR)$ is closed too.

  Consider the locally free sheaf $\kH=\bop_{k=m+1}^n \hom(b_k\kA,b_{k-1}\kR)$ and the projective
 space bundle $\mP(\kH)$ \cite[Section II.7]{ha}.  Every point $h\in\mP(\kH)$ defines a set of homomorphisms
 $h_k:b_k\kA(x)\to b_{k-1}\kR(x)$ (up to a homothety), where $x$ is the image of $h$ in $X$, and the points $h$
 such that $h_kh_{k+1}=0$ form a closed subset $\mD\subseteq \mP(\kH)$. We denote by $\pi$ the restriction onto
 $\mD$ of the projection $\mP(\kH)\to X$; it is a projective, hence closed mapping. Moreover, for every point $x\in X$
 the fibre $\pi^{-1}(x)$ is isomorphic to $\mD(\fB,\kA(x),\kR(x))$. Consider also the group variety $\kG$ over $X$:
 $\kG=\prod_{k=m}^n\mathrm{GL}_{b_k}(\kA)$. There is a natural action of $\kG$ on $\mD$ over $X$, and 
 the sets $\mD_i=\setsuch{z\in\mD}{\dim\kG z\le i}$ are closed in $\mD$. Therefore the sets $Z_i=\pi(\mD_i)$ are 
 closed in $X$, as well as $Z_{ij}=\setsuch{x\in Z_i}{\dim\pi^{-1}(x)\ge i+j}$. But $X_j(\kR)=\bigcup_iZ_{ij}$,
 thus it is also a closed set.
 \end{proof} 

 \section{Derived tame and wild algebras}
 \label{s4}

 To define derived tame and wild algebras we need consider families of complexes based on non-commutative algebras.
 As before, we assume that the field $\Mk$ is algebraically closed, though the definitions do not use this restriction.

  \begin{defin}\label{41}
  Let $\bA$ be a finite dimensional algebra over the field $\Mk$ with radical $\bR$ and $\bB$ be any $\Mk$-algebra.
 \begin{enumerate}
 \item
  A \emph{family of minimal $\bA$-complexes} based on $\bB$ is defined as a complex
  \begin{equation}\label{e41}
   \begin{CD}
   \kP\bul:\quad \kP_n@>{d_n}>>\kP_{n-1}@>{d_{n-1}}>>\dots\larr \kP_m ,  
 \end{CD} 
 \end{equation} 
 of finitely generated projective $\bB^\circ\*\bA$-modules such that $\im d_k\subseteq \kP_{k-1}\bR$ for all $k$.
 \item
  For a family \eqref{e41} and a finite dimensional (over $\Mk$) left $\bB$-module $L$ we denote by $\kP\bul(L)$
 the complex  
 $$ 
     \begin{CD}
   L\*_\kB\kP_n@>{1\*d_n}>>L\*_\bB\kP_{n-1}@>{1\*d_{n-1}}>>\dots\larr L\*_\bB\kP_m .
 \end{CD} 
 $$ 
 \item
  We call a family \eqref{e41} \emph{strict} if
    \begin{enumerate}
  \item 
    $\kP\bul(L)\simeq\kP\bul(L')$ \iff $L\simeq L'$;
  \item
    $\kP(L)$ is indecomposable \iff $L$ is indecomposable.
   \end{enumerate}
 \item
  We call $\bA$ \emph{derived wild} if for every finitely generated $\Mk$-algebra $\bB$ there is a strict family of
 minimal $\bA$-complexes based on $\bB$.
 \item
  We call $\bA$ \emph{derived tame} if there is a set $\dM$ of families of minimal $\bA$-complexes with the following
 properties:
  \begin{enumerate}
  \item 
   Every $\kP\bul\in\dM$ is based on a \emph{rational algebra} $\bB$, which means that $\bB\simeq\Mk[x,f(x)^{-1}]$ 
 for a non-zero polynomial $f(x)$. We define $\fR\bul(\kP\bul)$ as $\fR\bul(\kP\bul(L))$ for some (hence any) 
 one-dimensional $\bB$-module $L$.
  \item
   The set 
  $$
   \dM(\fR\bul)=\setsuch{\kP\bul\in\dM}{\fR\bul(\kP\bul)=\fR\bul}
  $$
  is finite for each $\fR\bul$.
   \item
    For every $\fR\bul$ all indecomposable minimal $\bA$-complexes of vector rank $\fR\bul$, except maybe finitely
 many isomorphism classes of such complexes, are isomorphic to $\kP\bul(L)$ for some $\kP\bul\in\dM$ and some
 $\bB$-module $L$.
  \end{enumerate}
 \end{enumerate}
 \end{defin} 

 \begin{erem}
  These definitions do not coincide but are easily seen to be equivalent to other used definitions of derived tame
 and wild algebras, for instance those of \cite{ge,br,gk}. As usually, to show that $\bA$ is derived wild it suffices to
 construct a strict family over one of specimen algebras such as free algebra $\Mk\langle x,y\rangle$, or polynomial
 algebra $\Mk[x,y]$, or power series algebra $\Mk[[x,y]]$.
 \end{erem} 

 The following proposition follows from elementary geometrical consideration like in \cite{d0}.

 \begin{prop}\label{42}
  \begin{enumerate}
\item 
  If  $\bA$ is derived tame, then $\pam(\fB,\bA)\le |\fB|\dim\bA$ for each sequence
 $\fB=(b_n,b_{n-1},\dots,b_m)$, where $|\fB|=\sum_{k=m}^n b_k$.
 \item
  If $\bA$ is derived wild, then there is a sequence $\fB$ such that $\pam(c\fB,\bA)\ge
 c^2$ for every integer $c$.
\end{enumerate}
 In particular, no algebra can simultaneously be both derived tame and derived wild.
 \end{prop} 

 In what follows we use the following supposition, which is believed by most experts.

 \begin{sps}\label{43}
  Every finite dimensional algebra is either derived tame or derived wild.%
\footnote{Now V.\,Bekkert and the author are preparing an article with a proof of this conjecture.}
 \end{sps} 

  \begin{corol}\label{44}
  Let $\kA$ be a family of algebras based on a algebraic variety $X$. Then $W=\setsuch{x\in X}{\kA(x)
 \text{ \em is derived wild}}$ is a union of a countable sequence of closed subsets. 
 \end{corol} 
 \begin{proof}
  Indeed $W=\bigcup_{c,\fB}W_{c,\fB}$, where $W_{c,\fB}=\setsuch{x\in X}{\pam(c\fB,\kA(x))>c|\fB|\rk\kA)}$.
 By Theorem \ref{31} all these subsets are closed in $X$. 
 \end{proof} 
 
  \begin{conj}\label{45}
 In the situation of Corollary \ref{44}, the set $W$ is always closed in $X$, or, the same, the set  
 $\setsuch{x\in X}{\kA(x) \text{ \em is tame}}$ is open.
 \end{conj} 

  \begin{corol}\label{46}
  Suppose that an algebra $\bA$, which is derived wild, \emph{degenerates} to another algebra $\ol\bA$,
 i.e. there is a family of algebras $\kA$ based on a variety $X$ such that $\kA(x)\simeq\bA$ for all $x$ from
 a dense open subset $U\sb X$ and there is a point $y\in X$ such that $\kA(y)\simeq\ol\bA$.
 Then $\ol\bA$ is also derived wild. 
 \end{corol} 
 
 \begin{rmk}\label{47}
 \def\bp{{\bullet}}
  Recently Th.\,Br\"ustle has announced a counter-example to the semi-continuity for derived categories
 (cf. \cite[Section 8.1]{br}). Namely, he claims that the derived wild algebra $\bA$ given by the quiver with relations
 $$ 
    
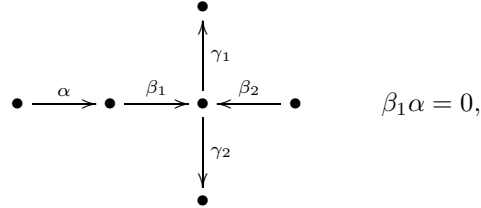
\xymatrix{
   && \bullet \\
   \bp \ar[r]^{\al} &\bp \ar[r]^{\be_1} &\bp \ar[u]_{\ga_1} \ar[d]^{\ga_2}
   & \bp \ar[l]_{\be_2} & {\be_1\al=0,} \\
   && \bullet
  }
 $$ 
 degenerates to the derived tame algebra $\ol\bA$ given by the quiver with relations 
 $$ 
   
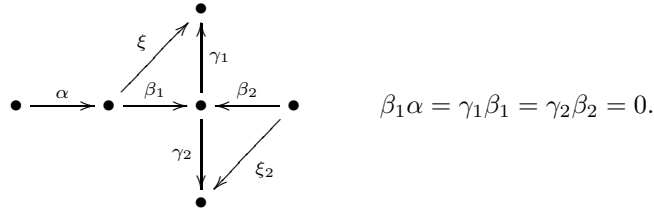
\xymatrix{
  & & \bp \\
  \bp \ar[r]^{\al} &\bp \ar[ur]^{\xi} \ar[r]^{\be_1} & \bp\ar[u]_{\ga_1} \ar[d]_{\ga_2}
  & \bp \ar[l]_{\be_2} \ar[dl]^{\xi_2} & {\be_1\al=\ga_1\be_1=\ga_2\be_2=0.}\\
 & & \bp
  } 
 $$ 
 As a matter of fact, it is not an example, since $\dim\bA=15$, while $\dim\ol\bA=16$,
 so the latter cannot be a degeneration of the former.
 \end{rmk}

\end{document}